\documentclass[a4paper]{article}
\usepackage{amsmath,amsthm,amssymb}
\usepackage[all]{xy} \SelectTips{cm}{}
\usepackage{rotating}
\theoremstyle{plain}      
    \newtheorem{theorem}{Theorem}[section]

\theoremstyle{definition}
    
    \newtheorem{definition-theorem}[theorem]{Definition-Theorem}

\theoremstyle{remark}
    \newtheorem*{remark}{Remark}
    \newtheorem*{note}{Note}
\newcommand{\A}{\ensuremath{\mathcal{A}}}
\newcommand{\B}{\ensuremath{\mathcal{B}}}
\newcommand{\C}{\ensuremath{\mathcal{C}}}
\newcommand{\D}{\ensuremath{\mathcal{D}}}
\newcommand{\E}{\ensuremath{\mathcal{E}}}
\newcommand{\X}{\ensuremath{\mathcal{X}}}
\newcommand{\V}{\ensuremath{\mathcal{V}}}
\newcommand{\Vect}{\ensuremath{\mathbf{Vect}}}

\newcommand{\End}{\ensuremath{\mathrm{End}}}
\newcommand{\fd}{\ensuremath{\mathrm{fd}}}

\DeclareMathOperator{\colim}{colim}

\newcommand{\dint}{\ensuremath{\displaystyle \int}}
\newcommand{\op}{\ensuremath{\mathrm{op}}}
\newcommand{\ox}{\ensuremath{\otimes}}

\newcommand{\<}{\ensuremath{\langle}}
\renewcommand{\>}{\ensuremath{\rangle}}

\newcommand{\ra}{\ensuremath{\rightarrow}}
\newcommand{\xra}{\ensuremath{\xrightarrow}}
\begin{document}
\title{On endomorphism algebras of functors with non-compact domain}
\author{Brian Day}
\maketitle
\begin{abstract}
As a development of~\cite{2} and~\cite{3}, we construct a ``VN-core'' in
$\Vect_k$ for each $k$-linear split-semigroupal functor from a suitable
monoidal category $\C$ to $\Vect_k$. The main aim here is to avoid the
customary compactness assumption on the set of generators of the domain
category $\C$ (cf.~\cite{3}).
\end{abstract}

%=========================================================================%
\section{Introduction}
%=========================================================================%

We propose the construction of a VN-core associated to each ($k$-linear)
split semigroupal functor $U$ from a suitable monoidal category $\C$ to
$\Vect_k$, where all our categories, functors, and natural transformations
are assumed to be $k$-linear, for a fixed field $k$. Essentially, the
category $\C$ must be equipped with a small ``$U$-generator'' $\A$ carrying
some extra duality information and with $UA$ still being finite dimensional
for all $A$ in $\A$.

We shall use the term ``VN-core'' (in $\Vect_k$) to mean a (usual)
$k$-semibialgebra $E$ together with a $k$-linear endomorphism $S$ such that
\[
    \mu(\mu \ox 1)(1 \ox S \ox 1)(1 \ox \delta)\delta = 1 : E \ra E.
\]
The VN-core is called ``antipodal'' if $S(xy) = Sy Sx$ (and $S(1) = 1$)
for all $x,y \in E$. This minimal type of structure is introduced here in
order to avoid compactness assumptions on the generator $\A \subset \C$ and,
at the same time, retain the ``fusion'' operator, namely
\[
    (\mu \ox 1)(1 \ox \delta) : E \ox E \ra E \ox E,
\]
satisfying the usual fusion equation~\cite{7}. Note that here the fusion
operator always has a partial inverse (see~\cite{1}).

In \S 2 we establish sufficient conditions on a functor $U$ in order that
\[
    \End^\vee U = \int^A (UA)^* \ox UA
\]
be a VN-core in $\Vect_k$ (following~\cite{2}). This core can be completed to a
VN-core $\End^\vee U \oplus k$ with a unit element. In \S 3 we give several
examples of suitable functors $U$ for the theory.

%=========================================================================%
\section{The construction of $\End^\vee U$}
%=========================================================================%

Let $\C = (\C,\ox,I)$ be a monoidal category and let
\[
    U : \C \ra \Vect
\]
be a functor with both a semigroupal structure, denoted
\[
    r = r_{C,D} : UC \ox UD \ra U(C \ox D),
\]
and a cosemigroupal structure, denoted
\[
    i = i_{C,D} : U(C \ox D) \ra UC \ox UD,
\]
such that $ri = 1$.

We shall suppose also that there exists a small full subcategory $\A$ of $\C$
with the properties:

\begin{enumerate}
\item $UA$ is finite dimensional for all $A \in \A$, 

\item $U$-density; the canonical map
\[
    \alpha_C : \int^A \C(A,C) \ox UA \ra UC
\]
is an isomorphism for all $C \in \C$,

\item there is an ``antipode'' functor
\[
    (-)^* : \A^\op \ra \A
\]
with a (``canonical'') map $e_A : A \ox A^* \ox A \ra A$ in $\C$ for each $A
\in \A$,

\item there is a natural isomorphism
\[
    u = u_A : U(A^*) \xra{\cong} U(A)^*,
\]

\item the following diagrams defining $\tilde{\tau}$, $\tilde{\rho}$ both
commute
\[
    \xygraph{{UA \ox U(A)^* \ox UA}="l"
        [r(2)u(1.2)]{UA \ox U(A^*) \ox UA}="t"
        [r(2)d(1.2)]{U(A \ox A^* \ox A)}="r"
        "t"[d(2.4)]{UA}="b"
        "l":@{.>}"r" ^-{\tilde{\tau}}
        "l":"t" ^-{1 \ox u^{-1} \ox 1}
        "t":"r" ^-{r_3}
        "l":"b" _-{e_{UA}}
        "r":"b" ^-{Ue_{A}}}
\]
and
\[
    \xygraph{{UA \ox U(A)^* \ox UA}="l"
        [r(2)u(1.2)]{UA \ox U(A^*) \ox UA}="t"
        [r(2)d(1.2)]{U(A \ox A^* \ox A)}="r"
        "t"[d(2.4)]{UA}="b"
        "r":@{.>}"l" _-{\tilde{\rho}}
        "t":"l" _-{1 \ox u \ox 1}
        "r":"t" _-{i_3}
        "l":"b" _-{e_{UA}}
        "r":"b" ^-{Ue_{A}}}
\]
where $e_{UA} = 1 \ox \mathrm{ev}$ in $\Vect$, and $r_3 i_3 = 1$.
\end{enumerate}

We now define the semibialgebra structure $(\End^\vee U,\mu,\delta)$ on 
\[
    \End^\vee U = \int^A U(A)^* \ox UA
\]
as in~\cite{2} \S 2, with the isomorphism of $k$-linear spaces
\[
    S = \sigma : \End^\vee U \ra \End^\vee U
\]
given (as in~\cite{2} \S 3) by the usual components
\[
    \xygraph{{U(A)^* \ox UA}="1"
        [r(3.5)]{U(A^*)^* \ox U(A^*)}="2"
        "1"[d]{U(A)^* \ox U(A)^{**}}="3"
        "2"[d]{U(A^*) \ox U(A^*)^*}="4"
        "1":"2" ^-{\sigma_A}
        "3":"4" ^-{u^{-1} \ox u^*}
        "1":"3" _-{1 \ox d}
        "4":"2" _-{c}}
\]
where $d$ is the canonical map from a vector space to its double dual.
Furthermore, each map
\[
    e_{UA} = 1 \ox \mathrm{ev}: UA \ox UA^* \ox UA \ra UA
\]
satisfies both the conditions
\begin{equation}\tag{E1}
    \vcenter{\xygraph{{UA}="1"
        [ru]{UA \ox UA^* \ox UA}="2"
        [rd]{UA}="3"
        "1":"2" ^-{n \ox 1}
        "2":"3" ^-{e_{UA}}
        "1":"3" _-1}}
\end{equation}
commutes, and
\begin{equation}\tag{E2}
    \vcenter{\xygraph{{UA^*}="1"
        [r(1.5)u]{UA^* \ox UA \ox UA^*}="2"
        [r(1.5)d]{UA^* \ox UA^{**} \ox UA^*}="3"
        "1":"2" ^-{1 \ox n}
        "2":"3" ^-{1 \ox d \ox 1}
        "1":"3" _-{e^*_{UA}}}}
\end{equation}
commutes, where $n = \mathrm{coev}:1 \ra UA \ox UA^*$ in $\Vect$.

Then we obtain:

\begin{theorem}
The structure $(\End^\vee U, \mu, \delta, S)$ is a VN-core in $\Vect_k$ which
can be completed to the VN-core $(\End^\vee U) \oplus k$.
\end{theorem}

\begin{proof}
The von Neumann axiom
\[
    \mu_3(1 \ox S \ox 1)\delta_3 = 1
\]
becomes the diagram (in which we have omitted ``$\ox$''):
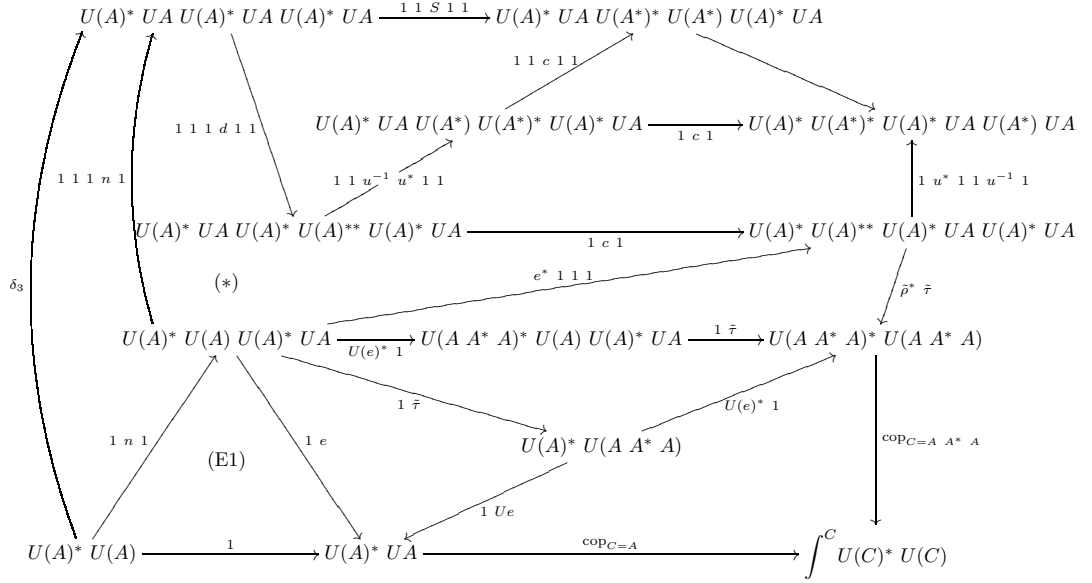
\[
    \scalebox{0.75}{
    \xygraph{{U(A)^*~UA~U(A)^*~UA~U(A)^*~UA}="1"
        [r(6)]{U(A)^*~UA~U(A^*)^*~U(A^*)~U(A)^*~UA}="2"
        "1"[d(3)r]{U(A)^*~UA~U(A)^*~U(A)^{**}~U(A)^*~UA}="3"
        [r(2.5)u(1.5)]{U(A)^*~UA~U(A^*)~U(A^*)^*~U(A)^*~UA}="4"
        [r(6)]{U(A)^*~U(A^*)^*~U(A)^*~UA~U(A^*)~UA}="5"
        "5"[d(1.5)]{U(A)^*~U(A)^{**}~U(A)^*~UA~U(A)^*~UA}="6"
        "3"[l(1)d(1.5)]{U(A)^*~U(A)~U(A)^*~UA}="7"
        [r(4.5)]{U(A~A^*~A)^*~U(A)~U(A)^*~UA}="8"
        [r(4.5)]{U(A~A^*~A)^*~U(A~A^*~A)}="10"
        "7"[r(5.2)d(1.5)]{U(A)^*~U(A~A^*~A)}="9"
        "7"[d(3)l(2)]{U(A)^*~U(A)}="11"
        "11"[r(4)]{U(A)^*~UA}="12"
        "12"[r(7)]{{\displaystyle \int^C} U(C)^*~U(C)}="13"
	    "7"[u(0.8)]{(*)}
	    "7"[d(1.7)]{\text{(E1)}}
        "1":"3" _-{1~1~1~d~1~1}
        "1":"2" ^-{1~1~S~1~1}
        "3":"4" |-{1~1~u^{-1}~u^*~1~1}
        "4":"5" _-{1~c~1}
        "4":"2" ^-{1~1~c~1~1}
        "2":"5" ^-{}
        "7":"8" _-{U(e)^*~1}
        "6":"5" _-{1~u^*~1~1~u^{-1}~1}
        "3":"6" _-{1~c~1}
        "7":@<8ex>@/^3ex/"1" ^-{1~1~1~n~1}
        "8":"10" ^-{1~\tilde{\tau}}
        "7":"9" _-{1~\tilde{\tau}}
        "9":"10" _-{U(e)^*~1}
        "9":"12" ^-{1~Ue}
        "11":"7" ^-{1~n~1}
        "11":"12" ^-{1}
        "7":"12" ^-{1~e}
        "6":"10" ^-{\tilde{\rho}^*~\tilde{\tau}}
        "10":"13" ^-{\textrm{cop}_{C = A~A^*~A}}
        "12":"13" ^-{\textrm{cop}_{C = A}}
        "11":@/^6ex/[u(7.3)] ^-{\delta_3}
	    "7":"6" ^-{e^*~1~1~1}
        }}
\]
where ($*$) is the exterior of the diagram
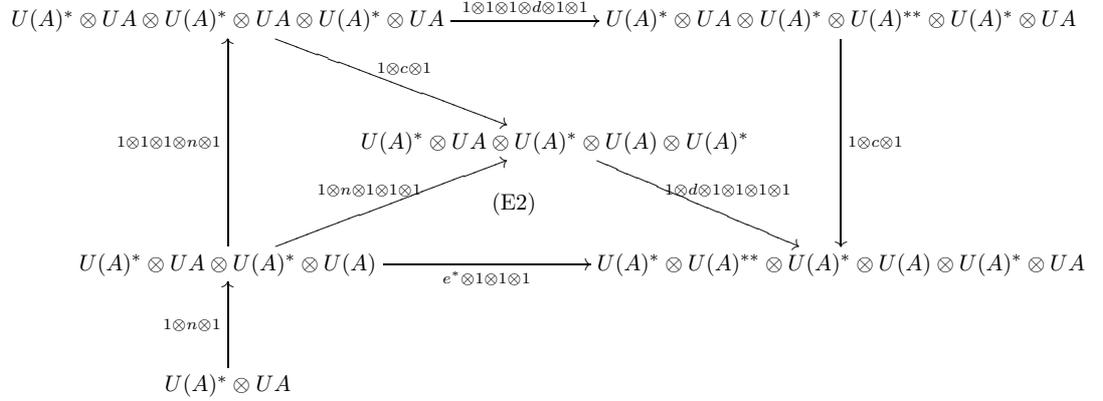
\[
    \scalebox{0.85}{
    \xygraph{{U(A)^* \ox  UA \ox U(A)^* \ox UA \ox U(A)^* \ox UA}="1"
        [r(7.5)]{U(A)^* \ox  UA \ox U(A)^* \ox U(A)^{**}  \ox U(A)^* \ox UA}="2"
        "1"[r(4)d(1.5)]{U(A)^* \ox  UA \ox U(A)^* \ox U(A) \ox U(A)^*}="3"
        "1"[d(3)]{U(A)^* \ox  UA \ox U(A)^* \ox U(A)}="4"
        "2"[d(3)]{U(A)^* \ox  U(A)^{**} \ox U(A)^* \ox U(A) \ox U(A)^* \ox UA}="5"
        "4"[d(1.5)]{U(A)^* \ox  UA}="6"
	"4"[r(3.5)u(0.75)]{\text{(E2)}}
	"1":"2" ^-{1 \ox 1 \ox 1 \ox d \ox 1 \ox 1}
	"1":"3" ^-{1 \ox c \ox 1}
	"2":"5" ^-{1 \ox c \ox 1}
	"3":"5" ^-{1 \ox d \ox 1 \ox 1 \ox 1 \ox 1}
	"4":"3" ^-{1 \ox n \ox 1 \ox 1 \ox 1}
	"4":"5" _-{e^* \ox 1 \ox 1 \ox 1}
	"4":"1" ^-{1 \ox 1 \ox 1 \ox n \ox 1}
	"6":"4" ^-{1 \ox n \ox 1}
    }}
\]
which commutes using (E2) and commutativity of
\[
    \xygraph{{\cdot}="1"
        [u(0.5)r(0.8)]{\cdot}="2"
        [r]{\cdot}="3"
        "1"[d(0.5)r(0.8)]{\cdot}="4"
        [r]{\cdot}="5"
        "1":"2" ^-n
        "2":"3" ^-{1~1~n}
        "3":"5" ^-{c}
        "1":"4" _-n
        "4":"5" _-{n~1~1}}
\]
\end{proof}

%=========================================================================%
\section{Examples}
%=========================================================================%

%=========================================================================%
\subsection{Example}\label{exam3.1}
%=========================================================================%

The first type of example is derived from the idea of a (contravariant)
involution on a (small) comonoidal category $\D$. This includes the doubles
$\D = \B^\op + \B$ and $\D = \B^\op \ox \B$ with their respective ``switch''
maps (where $\B$ is a given small comonoidal $\Vect_k$-category), or any
small comonoidal and compact-monoidal $\Vect_k$-category $\D$ (such as the
category $\mathbf{Mat}_k$ of finite matrices over $k$) with the tensor duals
of objects now providing an antipode on the comonoidal aspect of the structure
rather than on the monoidal part, or any $*$-algebra structure on a given
$k$-bialgebra (e.g., a $C^*$-bialgebra) with the $*$-operation providing
the antipode.

In each case, an \emph{even} functor from $\D$ to $\Vect$ is defined to be a
($k$-linear) functor $F$ equipped with a (chosen) dinatural isomorphism
\[
    F(D^*) \cong F(D).
\]
If we take the morphisms of even functors to be all the natural transformations
between them then we obtain a category
\[
    \E = \E(\D,\Vect).
\]
Let $\A = \E(\D,\Vect_\fd)_\mathrm{fs}$ be the full subcategory of $\E$
consisting of the finitely valued functors of finite support. While this
category is generally not compact, it has on it a natural antipode derived from
those on $\D$ and $\Vect_\fd$, namely
\[
    A^*(D) := A(D^*)^*.
\]
Of course, there are also examples where $\A$ is actually compact, such as
those where $\D$ is a Hopf algebroid, in the sense of~\cite{4}, with antipode
$(-)^* = S$, in which case each $A$ from $\D$ to $\Vect$ has a symmetry
structure on it.

Now let $\C$ be the full subcategory of $\E$ consisting of the small coproducts
in $\E$ of objects from $\A$. This category $\C$ is easily seen to be monoidal
under the pointwise convolution structure from $\D$, and the inclusion $\A
\subset \C$ is $U$-dense for the functor
\[
    U : \C \ra \Vect_k
\]
given by
\[
    U(C) = \sum_D C(D)
\]
which is split semigroupal with $UA$ finite dimensional for all $A \in \A$.
Moreover,
\begin{align*}
    U(A^*) &= \bigoplus_D A^*(D) \\
           &= \bigoplus_D A(D)^* \\
           &= U(A)^*,
\end{align*}
for all $A \in \A$. The conditions of (5) are easily verified if we define maps
\[
    e : A \ox A^* \ox A \ra A
\]
by commutativity of the diagrams
\[
    \xygraph{{A(D) \ox A^*(D) \ox A(D)}="1"
        [r(3)]{A(D)}="2"
        "1"[d]{A(D) \ox A(D)^* \ox A(D),}="3"
        "1":"2" ^-{e_D}
        "1":"3" _-{\cong}
        "3":"2" _-{1 \ox \mathrm{ev}}}
\]
where the exterior of
\[
    \xygraph{{A^*(D) \ox A(D)}="1"
        [d]{A(D)^* \ox A(D)}="2"
        [l(3)d]{A^*(E) \ox A(D)}="3"
	"2"[d]{A(E)^* \ox A(D)}="4"
	"2"[r(2.8)d]{k}="5"
        "4"[d]{A(E)^* \ox A(E)}="6"
        [d]{A^*(E) \ox A(E)}="7"
        "1":"2" _-{\cong}
        "1":"3" _-{A^*(f) \ox 1}
        "1":"5" ^-{\hat{e}}
        "3":"4" ^-{\cong}
        "4":"2" ^-{A(f)^* \ox 1}
        "2":"5" _-{\mathrm{ev}}
        "4":"6" _-{1 \ox A(f)}
        "6":"5" ^-{\mathrm{ev}}
        "6":"7" _-{\cong}
        "3":"7" _-{1 \ox A(f)}
        "7":"5" _-{\hat{e}}}
\]
commutes for all maps $f:D \ra E$ in $\D$ so that
\[
    e = 1 \ox \hat{e} : A \ox A^* \ox A \ra A \ox k \cong A
\]
is a genuine map in $\C$ when $\C$ is given the pointwise monoidal structure
from $\D$.
This completes the details of the general example.

%=========================================================================%
\subsection{Example}
%=========================================================================%

In the case where $k = \mathbb{C}$ and $\D$ has just one object $D$
whose endomorphism algebra is a $C^*$-bialgebra, we have a one-object
comonoidal category $\D$ with a $\mathbb{C}$-conjugate-linear antipode given by
the $*$-operation. Then the convolution $[\D, \mathbf{Hilb}_\fd]$, where
\[
    [\D, \mathbf{Hilb}_\fd] \subset [\D, \Vect_\mathbb{C}],
\]
is a monoidal category, with a $\mathbb{C}$-linear antipode given by
\[
    F^*(D) = F(D^*)^\circ
\]
where $H^\circ$ denotes the conjugate-transpose of $H \in \mathbf{Hilb}_\fd$.
We now interpret an even functor $F$ to be a functor equipped with a dinatural
isomorphism $F(D^*) \cong F(D)$ in $D \in \D$ which is
\emph{$\mathbb{C}$-linear}, so that $F^*(D) \cong F(D)^\circ$ for such a
functor.

Take $\A = \E(\D,\mathbf{Hilb}_\fd)$ and let $\C$ be the class of small
coproducts in $[\D,\Vect_\mathbb{C}]$ of the underlying
$[\D,\Vect_\mathbb{C}]$-representations of $A$'s in $\A$ (with the appropriate
maps). Each map
\[
    e : A \ox A^* \ox A \ra A
\]
in $\C$ is defined by the $\mathbb{C}$-linear components
\[
    e : A(D) \ox A^*(D) \ox A(D) \xra{1 \ox \hat{e}} A(D),
\]
where
\[
    \hat{e} : A^*(D) \ox A(D) \ra \mathbb{C}
\]
in $\Vect_\mathbb{C}$ comes from the $\mathbb{C}$-bilinear composite of two
maps which are both $\mathbb{C}$-linear in the first variable and
$\mathbb{C}$-linear in the second, namely
\[
    \xygraph{{A^*(D) \times A(D)}="1"
        [r(2)]{\mathbb{C}}="2"
        "1"[d]{A(D)^\circ \times A(D).}="3"
        "1":"2"
        "1":"3" _-{\cong}
        "3":"2" _-{\< -,- \>}}
\]
The remainder of this example is as seen before in Example~\ref{exam3.1}.

%=========================================================================%
\subsection{Example}
%=========================================================================%

Let $\V = (\V,\ox,I)$ be a (small) braided monoidal category and let $\B$ be
the $k$-linearization of $\mathbf{Semicoalg}(\V)$ with the monoidal structure
induced from that on $\V$.
By analogy with~\cite{5}, let $\X \subset \B$ be a finite full subcategory of
$\B$ with $I \in \X$ and $\X^\op$ promonoidal when
\begin{align*}
    p(x,y,z) &= \B(z, x \ox y) \\
    j(z) &= \B(z, I)
\end{align*}
for $x,y,z \in \X$.

For example (cf.~\cite{5}), one could take $\X$ to be a (finite) set of
non-isomorphic ``basic'' objects in some braided monoidal category $\V$, where
each $x \in \X$ has a coassociative diagonal map $\delta : x \ra x \ox x$.
However, we won't need the category $\X$ to be discrete or locally finite in
the following.

Now let $\C$ be the convolution $[\X^\op,\Vect]$ and let $\A =
[\X^\op,\Vect_\fd]$. The functor
\[
    U : \C \ra \Vect
\]
is defined by
\[
    U(C) = \bigoplus_x C(x),
\]
and the obvious inclusion $\A \subset \C$ is $U$-dense. If there is a canonical
(natural) retraction
\[
    \xygraph{{p(x,y,z) = \B(z, x \ox y)}="1"
        [r(4)]{\B(z,x) \ox \B(z,y),}="2"
	"1":@<-0.7ex>"2" _-{i_{x,y}}
	"2":@<-0.7ex>"1" _-{r_{x,y}}}
\]
derived from the semicoalgebra structures on $x,y,z$, then $U$ becomes a split
semigroupal functor via the structure maps
\[
    \xygraph{{U(C) \ox U(D)}="1"
        [r(5)]{U(C \ox D)}="2"
        "1"[d(1.3)]{\bigoplus_x C(x) \ox \bigoplus_y D(y)}="3"
        "2"[d(1.3)]{\bigoplus_z \dint^{xy} p(x,y,z) \ox C(x) \ox D(y)}="4"
        "3"[d(1.3)]{\bigoplus_z C(z) \ox D(z)}="5"
        "4"[d(1.3)]{\bigoplus_z \dint^{xy} \B(z,x) \ox \B(z,y) \ox C(x)\ox D(y),}="6"
        "1":@<3pt>"2" ^-r
        "2":@<3pt>"1" ^-i
        "1":@{=}"3"
        "2":@{=}"4"
        "5":@<3pt>"3" ^-{\Delta}
        "3":@<3pt>"5" ^-{\Delta^*}
        "6":@<3pt>"4" ^-{\text{``$r$''}}
        "4":@<3pt>"6" ^-{\text{``$i$''}}
        "6":"5" _-\cong}
\]
where the isomorphism follows from the Yoneda lemma, and $ri = 1$.

If $\X$ also has on it a duality
\[
    (-)^* : \X \ra \X^\op
\]
such that $x \cong x^{**}$, then, on defining
\[
    A^*(x) = A(x^*)^*,
\]
we obtain
\begin{align*}
    U(A^*) &= \bigoplus_x A^*(x) \\
           &= \bigoplus_x A(x^*)^* \\
           &\cong \bigoplus_x A(x)^* & \text{since $x \cong x^{**}$} \\
           &\cong U(A)^*,
\end{align*}
for $A \in \A$, in accordance with the fourth requirement on $U$.

Finally, to obtain a suitable map
\[
    e = 1 \ox \hat{e}: A \ox A^* \ox A \ra A \ox I \cong A,
\]
where $\hat{e}: A^* \ox A \ra I$, we suppose each $A$ in $\A$ has on it a
``dual coupling''
\[
    \chi = \chi_{xy} : A(x)^* \ox A(y) \ra \B(x^* \ox y,I).
\]
By considering the Yoneda expansion
\[
    A(x) \cong \int^z A(z) \ox \X(x,z)
\]
of the various functors $A$ in $\A = [\X^\op,\Vect_\fd]$, such a coupling
exists on each $A$ if we suppose merely that $\X$ itself is ``coupled'' by a
natural transformation
\[
    \chi : \X(y,z) \ra \X(x,z) \ox \B(x^* \ox y,I);
\]
or simply
\[
    \chi : \X(x,z)^* \ox \X(y,z) \ra \B(x^* \ox y,I),
\]
if $\X$ is locally finite. Then, the composite natural transformation
\[
    \xygraph{{A(x^*)^* \ox A(y) \ox \B(z, x \ox y)}="1"
        [d]{\B(x^{**} \ox y, I) \ox \B(z, x \ox y)}="2"
        [d]{\B(x \ox y, I) \ox \B(z, x \ox y)}="3"
        [d]{\B(z, I)}="4"
        "1":"2" ^-{\chi \ox 1}
        "2":"3" ^-\cong
        "3":"4" ^-{\text{comp'n}}}
\]
yields the map
\[
    \xygraph{{A^* \ox A}="1"
        [r(3.3)]{I}="2"
        "1"[d]{\dint^{xy} A^*(x) \ox A(y) \ox p(x,y,-)}="3"
        "2"[d]{\B(-,I)}="4"
        "1":"2" ^-{\hat{e}}
        "3":"4"
        "1":@{=}"3"
        "2":@{=}"4"}
\]
because $p(x,y,-) = \B(-,x \ox y)$ (by definition). Thus suitable conditions on
the coupling $\chi$ give (5).

\begin{remark}
Actually, this last example in which the basic promonoidal structure occurs as
a canonical retract of a comonoidal structure is typical of many other examples
which can be treated along similar lines.
\end{remark}

%=========================================================================%

%=========================================================================%

\bigskip
{\small\noindent
Department of Mathematics\\
Macquarie University\\
NSW, 2109, Australia}
\end{document}